\documentclass[12pt]{amsart}
\usepackage{latexsym}
\usepackage{amsmath, amssymb}
\newtheorem{theorem}{Theorem}[section]
\newtheorem{lms}{Lemma}[section]

\newtheorem{cors}{Corollary}[section]

\begin{document}
\title[Nuclear C$^*$-algebras]{On the classification problem \\ for nuclear C$^*$-algebras}

\author{Andrew S. Toms}
\date{January 31, 2005}
\address{Department of Mathematics and Statistics, York University, 4700 Keele St., Toronto, ON, Canada, M3J 1P3}
\email{atoms@mathstat.yorku.ca}

\begin{abstract}

We exhibit a counterexample to Elliott's classification conjecture
for simple, separable, and nuclear C$^*$-algebras whose construction is elementary,
 and demonstrate the necessity of extremely fine invariants
in distinguishing both approximate unitary equivalence 
classes of automorphisms of such algebras and isomorphism classes of the algebras themselves.
The consequences for the program to classify nuclear C$^*$-algebras
are far-reaching:  one has, among other things, that existing results on the classification
of simple, unital AH algebras via the Elliott invariant of $\mathrm{K}$-theoretic
data are the best possible, and that these cannot be improved by the addition
of continuous homotopy invariant functors to the Elliott invariant.

\end{abstract}

\maketitle

\section{Introduction}

Elliott's program to classify nuclear C$^*$-algebras via $\mathrm{K}$-theoretic
invariants (see \cite{El2} for an overview) has met with considerable success
since his seminal classification of approximately finite-dimensional (AF) algebras via their scaled, ordered
$\mathrm{K}_0$-groups (\cite{El1}).  Classification results 
of this nature are \emph{existence theorems} asserting that isomorphisms at the level of 
certain invariants for C$^*$-algebras in a class $\mathcal{B}$ are liftable to
$*$-isomorphisms at the level of the algebras themselves.  Obtaining such
theorems usually requires proving a  \emph{uniqueness theorem} 
for $\mathcal{B}$, i.e., a theorem which asserts that two $*$-isomorphisms between
members $A$ and $B$ of $\mathcal{B}$ which agree at the level of said invariants 
differ by a locally inner automorphism. 

Elliott's program began in earnest with his classification of simple circle algebras of real rank zero
in 1989 --- he conjectured shortly thereafter that the topological $\mathrm{K}$-groups, the
Choquet simplex of tracial states, and the natural connections between these objects 
would form a complete invariant for the class of separable, nuclear C$^*$-algebras.
This invariant came to be known simply as the Elliott invariant, denoted by $\mathrm{Ell}(\bullet)$.  
Elliott's conjecture held in the case of simple algebras throughout the 1990s, during which
time several spectacular classification results were obtained:  the Kirchberg-Phillips
classification of simple, separable, nuclear, and purely infinite (Kirchberg) C$^*$-algebras satisfying the 
Universal Coefficient Theorem, the Elliott-Gong-Li classification of simple unital AH
algebras of slow dimension growth, and Lin's classification of tracially AF algebras (see \cite{K}, \cite{EGL}, and
\cite{L}, respectively).
In 2002, however, R{\o}rdam constructed a simple, nuclear C$^*$-algebra containing both a finite
and an infinite projection (\cite{R1}).  Apart from answering negatively the question of whether 
simple, nuclear C$^*$-algebras have a type decomposition similar to that of
factors, his example provided the first counterexample to Elliott's conjecture
in the simple nuclear case;  it had the same Elliott invariant as a Kirchberg algebra
--- its tensor product with the Jiang-Su algebra $\mathcal{Z}$, to be precise --- yet was not purely infinite.  
It could, however, be distinguished from its Kirchberg twin by its (non-zero) real rank 
(\cite{R4}).  Later in the same year, the present author found independently a simple, nuclear, separable
and stably finite counterexample to Elliott's conjecture (\cite{T}).  This algebra could again be distinguished 
from its tensor product with the Jiang-Su algebra $\mathcal{Z}$ by its real rank.
These examples made it clear that the Elliott conjecture would not hold at its boldest, but 
the question of whether the addition of some small amount of new information to $\mathrm{Ell}(\bullet)$ 
could repair the defect in Elliott's conjecture remained unclear.  The counterexamples above suggested
the addition of the real rank, and such a modification would not have been without precedent:  
the discovery that the pairing between traces and the $\mathrm{K}_0$-group was necessary for
determining the isomorphism class of a nuclear C$^*$-algebra was unexpected, yet
the incorporation of this object into the Elliott invariant led to the classification of 
approximately interval (AI) algebras (\cite{E3}).  

The sequel clarifies the nature of the information not captured
by the Elliott invariant.  We exhibit a pair of simple, separable, nuclear, and non-isomorphic C$^*$-algebras
which agree  not only on $\mathrm{Ell}(\bullet)$, but also on a host of other invariants 
including the real rank and continuous (with respect to inductive sequences) homotopy
invariant functors.  The Cuntz semigroup, employed to distinguish our algebras, is thus
the minimum quantity by which the Elliott invariant must be enlarged in order
to obtain a complete invariant, but we shall see that the question of range for
this semigroup is out of reach.  Any classification result for C$^*$-algebras which
includes this semigroup as part of the invariant will therefore 
lack the impact of the Elliott program's successes --- 
the latter are always accompanied by range-of-invariant results.  
Our aim, however, is not to discourage work on the classification program.  It is to 
demonstrate unequivocally the need for a new regularity assumption in Elliott's 
program, as opposed to an expansion of the invariant.

Let $\mathcal{F}$ denote the following collection of invariants for C$^*$-algebras:
\begin{enumerate}
\item[$\bullet$] all homotopy invariant functors from the
category of C$^*$-algebras which commute with countable inductive limits;
\item[$\bullet$] the real rank (denoted by $\mathrm{rr}(\bullet)$);
\item[$\bullet$] the stable rank (denoted by $\mathrm{sr}(\bullet)$);
\item[$\bullet$] the Hausdorffized algebraic $\mathrm{K}_1$-group;
\item[$\bullet$] the Elliott invariant.
\end{enumerate}
Let $\mathcal{F}_{\mathbf{R}}$ be the subcollection of $\mathcal{F}$ obtained
by removing those continuous and homotopy invariant functors which do not have
ring modules as their target category.

Our main results are:
\begin{theorem}
There exists a simple, separable, unital, and nuclear C$^*$-algebra $A$ such that for any
UHF algebra $\mathcal{U}$ and any $F \in \mathcal{F}$ one has
\[
F(A) \cong F(A \otimes \mathcal{U}),
\]
yet $A$ and $A \otimes \mathcal{U}$ are not isomorphic.
$A$ is moreover an approximately homogeneous (AH) algebra, and
$A \otimes \mathcal{U}$ is an approximately interval (AI) algebra. 
\end{theorem}
\begin{theorem}
There exist a simple, separable, unital, and nuclear C$^*$-algebra $B$ and an automorphism $\alpha$
of $B$ of period two such that $\alpha$ induces the identity map on $F(B)$
for every $F \in \mathcal{F}_{\mathbf{R}}$, yet $\alpha$ is not locally inner.
\end{theorem}
Thus, both existence and uniqueness fail for simple, separable, and nuclear C$^*$-algebras
despite the scope of $\mathcal{F}$.

Recall that a C$^*$-algebra $A$ is said to be $\mathcal{Z}$-stable if it absorbs the Jiang-Su
algebra $\mathcal{Z}$ tensorially, i.e., $A \otimes \mathcal{Z} \cong A$.  ($\mathcal{Z}$-stability
is the regularity property alluded to above.)
Theorem 1.1, or rather, its proof, has two immediate corollaries which are of independent interest.
\begin{cors}
There exists a simple, separable, and nuclear C$^*$-algebra with unperforated
ordered $\mathrm{K}_0$-group whose Cuntz semigroup fails to be almost unperforated.
\end{cors}
\begin{cors}
Say that a simple, separable, nuclear, and stably finite C$^*$-algebra has property $(M)$ if it
has stable rank one, weakly unperforated topological K-groups, weak
divisibility, and property (SP).  Then, $(M)$ is strictly weaker than $\mathcal{Z}$-stability.
\end{cors}
Corollary 1.1 follows from the proof of Theorem 1.1, while Corollary 1.2 
follows from Corollary 1.1 and Theorem 4.5 of \cite{R3}.

 The counterexample to 
the Elliott conjecture constituted by Theorem 1.1 is more powerful and succinct than
those of \cite{R1} or \cite{T}:  $A$ and 
$A \otimes \mathcal{U}$ agree on the distinguishing invariant for the 
counterexamples of \cite{R1} and \cite{T} and a host of others including 
$\mathrm{K}$-theory with coefficients mod $p$, the homotopy groups of the
unitary group, the stable rank, and all $\sigma$-additive homologies and cohomologies from the
category of nuclear C$^*$-algebras (cf. \cite{Bl1});  $A$ and
$A \otimes \mathcal{U}$ are simple, unital inductive limits
of homogeneous algebras with contractible spectra, a class of algebras which forms
the weakest and most natural extension of the very slow dimension growth AH algebras classified in \cite{EGL};   
both $A$ and $A \otimes \mathcal{U}$ are stably finite, weakly divisible, and have 
property (SP), minimal stable rank, and next-to-minimal
real rank;  the proof of the theorem is elementary compared to 
the intricate constructions of \cite{R1} and \cite{T}, and demonstrates the necessity of
a distinguishing invariant for which no range results can be expected.  
Furthermore, one has in Theorem 1.2 a companion lack-of-uniqueness result.  
Together with Theorem 1.1, this yields what might be called a categorical
counterexample --- the structure of the category whose objects are isomorphism classes of
simple, separable, nuclear, stably finite C$^*$-algebras (let alone just
nuclear algebras) and whose morphisms are locally inner equivalence classes
of $*$-isomorphisms cannot be determined by $\mathcal{F}$.

The paper is organized as follows:  Section 2 fixes notation and reviews the definition
of the Cuntz semigroup $W(\bullet)$;  in Section 3 we prove
Theorem 1.1;  in Section 4 we prove Theorem 1.2;  Section 5 demonstrates
the complexity of the Cuntz semigroup, and discusses the relevance of 
$\mathcal{Z}$-stability to the classification program.

\emph{Acknowledgements}.  The author would like to thank Mikael R{\o}rdam 
both for suggesting the search for the automorphisms of Theorem 1.2 and for
several helpful discussions, S{\o}ren Eilers and Copenhagen University for
their hospitality in 2003, and George Elliott for his hospitality and comments 
at the Fields Institute in early 2004, where some of the work on Theorem 1.2 
was carried out.  This work was supported by an NSERC Postdoctoral Fellowship and
by a University of New Brunswick grant.

\section{Preliminaries}

For the remainder of the paper, let $\mathrm{M}_n$ denote
the $n \times n$ matrices with complex entries, and let $\mathrm{C}(X)$ denote
the continuous complex-valued functions on a topological space $X$.

Let $A$ be a C$^*$-algebra.  We recall the definition of the Cuntz semigroup $W(A)$ from \cite{C}.  (Our 
synopsis is essentially that of \cite{R3}.)  Let
$\mathrm{M}_n(A)^+$ denote the positive elements of $\mathrm{M}_n(A)$,
and let $\mathrm{M}_{\infty}(A)^+$ be the disjoint union $\cup_{i=n}^{\infty} 
\mathrm{M}_n(A)^+$.  For $a \in \mathrm{M}_n(A)^+$ and $b \in \mathrm{M}_m(A)^+$
set $a \oplus b = \mathrm{diag}(a,b) \in \mathrm{M}_{n+m}(A)^+$, and write 
$a \precsim b$ if there is a sequence $\{ x_k \}$ in $\mathrm{M}_{m,n}(A)$ such
that $x_k^* b x_k \rightarrow a$.  Write $a \sim b$ if $a \precsim b$ and 
$b \precsim a$.  Put $W(A) = \mathrm{M}_{\infty}(A)^+ / \sim$, and let
$\langle a \rangle$ be the equivalence class containing $a$.  Then,
$W(A)$ is a positive ordered abelian semigroup when equipped with the
relations:
\begin{displaymath}
\langle a \rangle + \langle b \rangle = \langle a \oplus b \rangle, \ \ \ \ \
\langle a \rangle \leq \langle b \rangle \Longleftrightarrow a \precsim b, \ \ \ \ \
a,b \in \mathrm{M}_{\infty}(A)^+.
\end{displaymath}
The relation $\precsim$ reduces to Murray-von Neumann comparison when $a$ and $b$
are projections.

We will have occasion to use the following simple lemma in the sequel:
\begin{lms}
Let $p$ and $q$ be projections in a C$^*$-algebra $D$ such that
\[
||xpx^*-q|| < 1/2
\]
for some $x \in D$.  Then, $q$ is equivalent to 
a subprojection of $p$.
\end{lms}

\begin{proof} We have that
\[
\sigma(xpx^*) \subseteq (-1/2,1/2) \cup (1/2,3/2),
\]
and that $\sigma(xpx^*)$ contains at least one point from $(1/2,3/2)$. 
The C$^*$-algebra generated by $xpx^*$ contains a non-zero projection, say $r$, 
represented (via the functional calculus) by the function $r(t)$ on $\sigma(xpx^*)$
which is zero when \mbox{$t \in (-1/2,1/2)$} and one otherwise.  This projection
is dominated by 
\[
2 xpx^* = \sqrt{2} x p x^* \sqrt{2}.
\]
By the functional
calculus one has $||xpx^* - r|| < 1/2$, so that $||r-q||<1$.  Thus,
$r$ and $q$ are Murray-von Neumann equivalent.  By the definition of 
Cuntz equivalence we have $\sqrt{2} x p x^* \sqrt{2} \precsim p$, so 
that $q \sim r \precsim p$ by transitivity.  Cuntz comparison
agrees with Murray-von Neumann comparison on projections, and the lemma
follows. \end{proof}

\section{The proof of Theorem 1.1}

\begin{proof}
We construct $A$ as an inductive 
limit $\lim_{i \to \infty}(A_i,\phi_i)$ where, for each $i \in \mathbb{N}$, $A_i$ is of the form 
\[
\mathrm{M}_{k_i} \otimes \mathrm{C}\left([0,1]^{6 (\Pi_{j \leq i} n_j)}\right), \ \ n_i, k_i \in
\mathbb{N}, 
\]
and $\phi_i$ is a unital
$*$-homomorphism.  Our construction is essentially that of \cite{V1}.  Put $k_1 = 4$, $n_1 = 1$, and $N_i = \Pi_{j \leq i} n_j$.
Let 
\[
\pi_l^i: [0,1]^{6 N_i} \to [0,1]^{6 N_{i-1}}, 
\ \ l \in \{1,\ldots,n_i\}, 
\]
be the co-ordinate projections, and let $f \in A_{i-1}$.  Define $\phi_{i-1}$ by
\[
\phi_{i-1}(f)(x) = \mathrm{diag}\left(f(\pi_1^i(x)),\ldots,f(\pi_{n_i}^i(x)),f(x_1^{i-1}),
\ldots,f(x_{m_i}^{i-1})\right),
\]
where $x_1^{i-1},\ldots,x_{m_i}^{i-1}$ are points in $X_{i-1} \stackrel{\mathrm{def}}{=} 
[0,1]^{6 N_{i-1}}$.  With $m_i = i$, the $x_1^{i-1},\ldots,x_{m_i}^{i-1}$, $i \in
\mathbb{N}$, can be chosen so as to make $\lim_{i \to \infty}
(A_i,\phi_i)$ simple (cf. \cite{V2}).  The multiplicity of $\phi_{i-1}$ is 
$n_i + m_i$ by construction.  We impose two conditions on the $n_i$ and $m_i$:  first, $n_i \gg m_i$
as $i \to \infty$, and second, given any natural number r, there is an $i_0 \in \mathbb{N}$
such that $r$ divides $n_{i_0} + m_{i_0}$.  

Note that $(\mathrm{K}_0 A_i,{\mathrm{K}_0^+ A_i},[1_{A_i}]) = (\mathbb{Z},\mathbb{Z}^+,k_i)$ since 
$X_i$ is contractible for all $i \in \mathbb{N}$.  The second condition on the $n_i$
above implies that
\[
(\mathrm{K}_0 A, \mathrm{K}_0 A^+, [1_{A}]) = \lim_{i \to \infty}(\mathrm{K}_0 A_i,{\mathrm{K}_0 A_i^+}, [1_{A_i}])
\cong (\mathbb{Q},\mathbb{Q}^+, 1).
\]
Since $\mathrm{K}_1 A_i = 0$, $i \in \mathbb{N}$, we have $\mathrm{K}_1 A = 0$.  Thus, 
$A$ has the same Elliott invariant as some AI algebra, say $B$.  Tensoring $A$ with
a UHF algebra $\mathfrak{U}$ does not disturb the $\mathrm{K}_0$-group or the tracial
simplex ($\mathfrak{U}$ has a unique normalized tracial state).  The tensor product
$A \otimes \mathfrak{U}$ is a simple, unital AH algebra with very slow dimension growth
in the sense of \cite{EGL},
and is thus isomorphic to $B$ by the classification theorem of \cite{EGL}.

Let us now prove that $A$ and $B$ are shape equivalent.  By the range-of-invariant
theorem of \cite{Th}
we may write $B$ as an inductive limit of full matrix algebras over the 
closed unit interval (as opposed to direct sums of such), say 
\[
B \cong \lim_{i \to \infty}(B_i,\psi_i).
\]
From $\mathrm{K}$-theory considerations we may assume that $B_i = \mathrm{M}_{k_i}
\otimes \mathrm{C}([0,1])$, i.e., that the dimension of the unit of $B_i$ is the 
same as the dimension of the unit of $A_i$.  Let $s_i = \mathrm{mult}\phi_i
= \mathrm{mult}\psi_i$.  Define maps
\[
\eta_i:A_i \to B_{i+1}, \ \ \eta_i(f) = \bigoplus_{j=1}^{s_i} f((0,\ldots,0))
\]
and
\[
\gamma_i:B_i \to A_i, \ \ \gamma_i(g) = g(0).
\]
Both $\gamma_{i+1} \circ \eta_i$ and $\eta_{i} \circ \gamma_{i-1}$
are diagonal maps, and so are homotopic to $\phi_i$ and $\psi_i$, respectively,
since $[0,1]$ and $X_i$ are contractible.

Finally, $A$ has stable rank one and real rank one by \cite{V2}, and therefore
so also does $B$.

To complete the proof of the theorem, we must show that $A$ and $B$
are non-isomorphic.  Since $B$ is approximately divisible, we have that $W(B)$ is 
almost unperforated, i.e., if $mx \precsim ny$ for natural numbers $m > n$ and elements
$x,y \in W(B)$, then $x \precsim y$ (\cite{R2}).  We claim that the Cuntz semigroup 
of $A$ fails to be almost unperforated.  We proceed by extending Villadsen's Euler class 
obstruction argument (cf. \cite{V1}, \cite{V2}) to positive elements of a particular form.  

To show that $W(A)$ fails to be almost unperforated, it will suffice to 
exhibit positive elements $x,y \in A_1$ such that, for all $i \in \mathbb{N}$, 
for some $\delta > 0$
\[
m \langle \phi_{1i}(x) \rangle \precsim n \langle \phi_{1i}(y) \rangle, \ \ m > n, \ \ m,n \in \mathbb{N}
\]
and
\[
||r \phi_{1i}(y) r^* - \phi_{1i}(x)|| > \delta, \ \ \forall r \in A_i, \ \ \forall i \in \mathbb{N}.
\]
The second statement is stronger than the requirement
that $\langle \phi_{1i}(x) \rangle$ is not less than $\langle \phi_{1i}(y) 
\rangle$ in $W(A_i)$, since $W(\bullet)$ does not commute with inductive
limits.  Clearly, we need only establish this second statement over
some closed subset $Y$ of the spectrum of $A_i$. 

If $a \in \mathrm{M}_n \otimes \mathrm{C}(X)$ is a constant positive element and
$X$ is compact, then $\langle a \rangle$ is the class of a projection in
$W(\mathrm{M}_n \otimes \mathrm{C}(X))$.  Indeed, $a$ is unitarily equivalent (hence Cuntz
equivalent) to a diagonal positive element:
\[
uau* = \mathrm{diag}(a_1,\ldots,a_m,0,\ldots,0), \ \ \mathrm{some} \ u \in \mathcal{U}(\mathrm{M}_n),
\]
where $a_l \neq 0$, $l \in \{1,\ldots,m\}$.  Let $r = \mathrm{diag}(a_1^{-1},\ldots,a_m^{-1},0,\ldots,0)$.
Then, 
\[
r^{1/2} u a u* r^{1/2} = (r^{1/2}u)a(r^{1/2}u)^* = \mathrm{diag}(\underbrace{1,\ldots,1}_{m \ \mathrm{times}},0,\ldots,0).
\]

Set 
\[
S \stackrel{\mathrm{def}}{=} \left\{\overline{x} \in [0,1]^3 \ : \ \frac{1}{8} < 
\mathrm{dist}\left(\overline{x},\left(\frac{1}{2},\frac{1}{2},\frac{1}{2}\right)\right) < \frac{3}{8}\right\}.
\]
Note that $\mathrm{M}_4(\mathrm{C}_0(S \times S))$ is a hereditary subalgebra
of $A_1$.
Let $\xi$ be a line bundle over $\mathrm{S}^2$ with non-zero Euler class 
(the Hopf line bundle, for instance).  Let $\theta_1$ denote
the trivial line bundle.  By Lemma 1 of \cite{V2}, we have that 
$\theta_1$ is not a sub-bundle of $\xi \times \xi$ over $\mathrm{S}^2 \times
\mathrm{S}^2$.  Both $\xi \times \xi$ and $\theta_1$ can be considered
as projections in $\mathrm{M}_4(\mathrm{S}^2 \times \mathrm{S}^2)$.  
By Lemma 2.1 we have
\[
||x (\xi \times \xi) x^* - \theta_1|| \geq 1/2, \ \forall x \in \mathrm{M}_4(\mathrm{S}^2 \times \mathrm{S}^2).
\]
On the other hand, the stability properties of vector bundles imply that
\[
11 \langle \theta_1 \rangle \leq 10 \langle \xi \times \xi \rangle.
\]

Consider the closure $S^-$ of $S \subseteq [0,1]^3$, and let $\tau$ be
the projection of $S^-$ onto 
\[
S_{1/4} \stackrel{\mathrm{def}}{=} \left\{\overline{x} \in S \ : \ \mathrm{dist}\left(\overline{x},
\left(\frac{1}{2},\frac{1}{2},\frac{1}{2}\right)\right)= \frac{1}{4}\right\} \subseteq 
[0,1]^3
\]
along rays emanating from $(1/2,1/2,1/2) \in [0,1]^3$.
Let $\tau^*(\xi)$ be the pullback of $\xi$ via $\tau$.
Define a function $f \in \mathrm{C}_0(S \times S)$ by
\[
f(\overline{x}) = 8 \mathrm{dist}(\overline{x},S_{1/4}).
\]
Note that $f$ takes the value $1$ on $S_{1/4}$.  By Lemma 2.1 we have
\[
||xf (\tau^*(\xi) \times \tau^*(\xi))x^* - f \theta_1|| \geq 1/2 
\]
for any $x \in A_1$ --- one simply restricts to $S_{1/4} \times S_{1/4} \subseteq S \times S$.
We may pull the inequality 
\[
11 \langle \theta_1 \rangle \leq 10 \langle \xi \times \xi \rangle.
\]
back via $\tau$ to conclude that
\[
11 \langle \theta_1 \rangle \leq 10 \langle \tau^*(\xi) \times \tau^*(\xi) \rangle.
\]
This last inequality is equivalent to the existence of a sequence $(r_j)$ in
the appropriately sized matrix algebra over $\mathrm{C}(S^- \times S^-)$ with the
property that 
\[
r_j \left( \oplus_{i=1}^{10} \tau^*(\xi) \times \tau^*(\xi) \right) r_j^* \stackrel{j \to \infty}{\longrightarrow} \theta_{11}.
\]
Since $f$ is central in $\mathrm{C}_0(S \times S)$, we have that
\[
r_j \left( \oplus_{i=1}^{10} f (\tau^*(\xi) \times \tau^*(\xi)) \right) r_j^* \stackrel{j \to \infty}{\longrightarrow} f  \theta_{11}.
\]
In other words,
\[
11 \langle f \theta_1 \rangle \leq 10 \langle f (\tau^*(\xi) \times \tau^*(\xi)) \rangle
\]
and $W(A_1)$ fails to be weakly unperforated.

Since
\[
11 \langle \phi_{1i}(f \theta_1) \rangle \leq 10 \langle \phi_{1i}(f (\tau^*(\xi) \times \tau^*(\xi))) \rangle
\]
via $\phi_{1i}(r_j)$, we need only show that 
\[
\left||x \phi_{1i}(f (\tau^*(\xi) \times \tau^*(\xi))) x^* - \phi_{1i}(f \theta_1) \right|| \geq 1/2 
\]
for each natural number $i$ and any $x \in A_i$.  Fix $i$.
One can easily verify that the restriction of $\phi_{1i}(f \cdot \tau^*(\xi) \times \tau^*(\xi))$
to $(S^-)^{2 N_i} \subseteq [0,1]^{6 N_i}$ is
\[
(\tau^*(\xi) \times \tau^*(\xi))^{\times N_i} \oplus f_{\theta_{l}},
\]
where $f_{\theta_l}$ is a constant positive element of rank $l$ (hence Cuntz equivalent 
to $\theta_l$), and the direct sum decomposition separates the summands of $\phi_{i-1}$
which are point evaluations from those which are not.
The similar restricted decomposition of $\phi_{1i}(f \cdot \theta_1)$ is
\[
\theta_{k-l/2} \oplus g_{\theta_{l/2}},
\]
where $g_{\theta_{l/2}}$ is a constant positive element Cuntz equivalent to a trivial projection of
dimension $l/2$, and $k$ is greater than $3l/2$ (this last inequality follows from the fact 
that $n_i \gg m_i$).
Suppose that there exists $x \in A_i|_{(S^-)^{2 N_i}}$ such that 
\[
||x ((\tau^*(\xi) \times \tau^*(\xi))^{\times N_i} \oplus f_{\theta_{l}}) x^* - 
\theta_{k-l/2} \oplus g_{\theta_{l/2}}|| < 1/2.
\]
Recall that 
\[
(\tau^*(\xi) \times \tau^*(\xi))^{\times N_i} \oplus f_{\theta_{l}} =
a((\tau^*(\xi) \times \tau^*(\xi))^{\times N_i} \oplus \theta_{l})a
\]
for some positive $a \in A_i$.   Cutting down by $\theta_{k-l/2}$, we have 
\[
||\theta_{k-l/2} x a((\tau^*(\xi) \times \tau^*(\xi))^{\times N_i} \oplus \theta_l)a x^* \theta_{k-l/2} - 
\theta_{k-l/2}|| < 1/2.
\]
By Lemma 2.1, we must conclude that 
\[
\theta_{k-l/2} \precsim (\tau^*(\xi) \times \tau^*(\xi))^{\times N_i} \oplus \theta_l
\]
over $(S^-)^{2 N_i}$.  
But this is impossible by Lemma 1 of \cite{V2}.  Hence
\[
||x(\phi_{1i}(f \cdot \tau^*(\xi) \times \tau^*(\xi)))x^* - \phi_{1i}(f \cdot \theta_1)|| \geq 1/2 \ \forall x \in A_i, 
\]
as desired. \end{proof}

\section{The proof of Theorem 1.2}

\begin{proof}
We perturb the construction of a simple, unital AH algebra by Villadsen (\cite{V1}) 
to obtain the algebra $B$ of Theorem 1.2, and construct $\alpha$ as an
inductive limit automorphism.
Let $X$ and $Y$ be compact connected Hausdorff spaces, and let
$\mathcal{K}$ denote the C$^*$-algebra of compact operators on
a separable Hilbert space.  Projections in the C$^*$-algebra $\mathrm{C}(Y) \otimes \mathcal{K}$
can be identified with finite-dimensional complex vector bundles over $Y$,
and two such bundles are stably isomorphic if and only if the
corresponding projections in $\mathrm{C}(Y) \otimes \mathcal{K}$ have the
same $\mathrm{K}_0$-class.

Given a set of mutually orthogonal projections 
\[
P= \{ p_1,\ldots,p_n\} \subseteq \mathrm{C}(Y) \otimes \mathcal{K}
\]
and continuous maps $\lambda_i: Y
\rightarrow X$, $1 \leq i \leq n$, one may define a \mbox{$*$-homomorphism}
\begin{displaymath}
\lambda: \mathrm{C}(X) \rightarrow \mathrm{C}(Y) \otimes
\mathcal{K}, \ \ f \rightarrow \bigoplus_{i=1}^{n} (f \circ \lambda_i)p_i.
\end{displaymath}
A \mbox{$*$-homomorphism} of this form is called \emph{diagonal}.  We say that $\lambda$
comes from the set ${\{(\lambda_i,p_i)\}}_{i=1}^{n}$.

Let $\mathrm{I}$ denote the closed unit interval in $\mathbb{R}$, and put
\begin{displaymath}
X_{i} = \mathrm{I} \times \mathrm{CP}^{\sigma(1)} \times
\mathrm{CP}^{\sigma(2)} \times \cdots \times \mathrm{CP}^{\sigma(i)},
\end{displaymath}
where the $\sigma(i)$ are natural numbers to be specified.
Let
\begin{displaymath}
\pi_{i+1}^1:X_{i+1} \to X_i; \ \
\pi_{i+1}^{2}:X_{i+1} \to \mathrm{CP}^{\sigma(i+1)}
\end{displaymath}
be the co-ordinate projections.
Let $B_{i} = p_{i}(\mathrm{C}(X_{i}) \otimes \mathcal{K})p_{i}$,
where $p_{i}$ is a projection in
$\mathrm{C}(X_{i}) \otimes \mathcal{K}$ to be specified.  The algebra
$B$ of Theorem 1.2 will be realized as the inductive limit of the $B_i$
with diagonal connecting \mbox{$*$-homomorphisms} $\gamma_i:B_i \to
B_{i+1}$.

Let $p_{1}$ be a projection corresponding to the vector bundle
\begin{displaymath}
\theta_{1} \times \xi_{\sigma(1)},
\end{displaymath}
over $X_{1}$, where $\theta_{1}$ denotes the trivial complex line bundle, $\xi_k$
denotes the universal line bundle over $\mathrm{CP}^{k}$ for a given
natural number $k$, and $\sigma(1) = 1$.  Put $\eta_i =
\pi_{i}^{2*}(\xi_{\sigma(i)})$.

We now
specify, inductively, the maps $\gamma_{i}:B_{i} \to B_{i+1}$.
Let $\tilde{\psi}$ be the homeomorphism of $\mathrm{I}$ given by
\[
\tilde{\psi}(x) = 1-x.
\]
Abusing notation, we will also take $\tilde{\psi}$ be the homeomorphism of $X_i \stackrel{\mathrm{def}}{=} \mathrm{I} \times Y_i$
given by $(x,y) \mapsto (\tilde{\psi}(x),y)$.  Choose a dense sequence 
$(z_i^l)_{l=1}^{\infty}$ in $X_i$ and choose for each $j= 1,2,\ldots,i+1$ a point $y_i^j \in X_i$ 
such that
$y_i^{i+1} = z_i^1$, $y_i^i = z_i^2$ and $\pi_{j+1}^1 \circ \pi_j^1 \circ 
\cdots
\circ \pi_i^1(y_i^j) = z_j^{i-j+2}$ for $1 \leq j \leq i-1$.  Let 
\begin{displaymath}
\tilde{\gamma_i}:
\mathrm{C}(X_i \otimes \mathcal{K}) \longrightarrow \mathrm{C}(X_{i+1} 
\otimes \mathcal{K})
\end{displaymath}
be a diagonal $*$-homomorphism coming from 
\begin{displaymath}
(\pi_{i+1}^1,\theta_1) \cup 
\{(y_i^j,\eta_{i+1})\}_{j=1}^{i+1} \cup \{(\tilde{\psi}(y_i^j),\eta_{i+1})\}_{j=1}^{i+1}.  
\end{displaymath}

Let 
$\tilde{\gamma_{1i}}$ be the composition
$\tilde{\gamma_i} \circ \cdots \circ \tilde{\gamma_1}$, and put $p_{i+1} 
= \tilde{\gamma_{1i}}(p_1)$
for all natural numbers $i$.  Let
$\gamma_i:B_i \to B_{i+1}$ be the restriction of 
$\tilde{\gamma_i}$.  Let
$B = \lim_{\rightarrow}(B_i,\gamma_i)$.  It follows from \cite{V1} that $B$ 
is simple, unital
AH-algebra.  (Apart from the choice of point evaluations in the 
$\tilde{\gamma_i}$,
the construction above is precisely that of~\cite{V1}.  The reason for the 
specific choice
of point evaluations will be made clear shortly.)

Straightforward calculation shows that the projection $p_i \in B_i$ 
corresponds to a complex vector bundle over $X_i$ of the
form $\theta_1 \oplus \omega_i$.  In fact,
with \mbox{$X_i = \mathrm{I} \times Y_i$} and with $\tau^i_1$, $\tau_2^i$ the 
co-ordinate projections,
we have that $\omega_i = \tau_2^{i*}(\tilde{\omega_i})$ for a vector bundle 
$\tilde{\omega_i}$ over $Y_i$.  
Thus, the homeomorphism $\tilde{\psi}$
of $X_i$ fixes $p_i$, and so induces and automorphism $\psi_i$ of $B_i$.  

Let $\pi_{im}^1$ be the composition $\pi_m^1 \circ \cdots \circ \pi_{i+1}^1$.
Let $f \in B_i$.  Then, with $(x,y)$ an element of $X_{i+1} = X_i \times 
\mathrm{CP}^{\sigma(i+1)}$,
we have
\begin{displaymath}
\gamma_i(f)(x,y) = f(\pi_{i+1}^1(x)) \oplus \left(\bigoplus_{j=1}^{i+1} 
f(\tilde{\psi}(y_i^j)) \otimes \eta_{i+1} \oplus f(y_i^j) \otimes \eta_{i+1} \right),
\end{displaymath}
so that
\begin{displaymath}
\psi_{i+1} \left(\gamma_i(f)(x,y)\right) = f\left(\tilde{\psi}(\pi_{i+1}^1(x))\right)
\oplus \left(\bigoplus_{j=1}^{i+1}  
f(\tilde{\psi}(y_i^j)) \otimes \eta_{i+1} \oplus f(y_i^j) \otimes \eta_{i+1} \right).
\end{displaymath}
On the other hand we have
\begin{displaymath}
\gamma_i \circ \psi_i(f)(x,y) = f\left(\tilde{\psi}(\pi_{i+1}^1(x))\right) \oplus 
\left(\bigoplus_{j=1}^{i+1}  
f(\tilde{\psi}(y_i^j)) \otimes \eta_{i+1} \oplus f(y_i^j) \otimes \eta_{i+1} \right).
\end{displaymath}
Thus, $\gamma_i \circ \psi_i$ and $\psi_{i+1} \circ \gamma_i$ 
differ only in the order of their direct summands, and so are
unitarily equivalent.  The unitary implementing this equivalence squares
to the identity.  Conjugating $\psi_{i+1}$ by said unitary element, we may assume that
$\gamma_i \circ \psi_{i} = \psi_{i+1} \circ \gamma_i$.  This process may
be repeated inductively for $\psi_m$, $m > i$, yielding an inductive limit
automorphism $\alpha$ of $B$ via the $\psi_i$.

We now show that $\alpha$ is not locally inner, yet induces
the identity map on $\mathrm{Inv}_F$ for any $F \in \mathcal{F}$.
Recall that the Euler class $e(\omega)$ of a complex vector bundle $\omega$ over a 
connected finite
CW-complex $X$ is an element of $\mathrm{H}^{2 \mathrm{dim} \omega}(X)$.  
For a trivial
complex vector bundle $\theta_l$ of dimension $l \in \mathbb{N}$ we 
have $e(\theta_l)=0$.
We also have $e(\omega_1 \oplus \omega_2) = e(\omega_1) \cdot 
e(\omega_2)$ for two
complex vector bundles $\omega_1$ and $\omega_2$ over $X$, where the product 
is the cup product in the integral cohomology ring $\mathrm{H}^*(X)$.  
Thus, if $e(\omega) \neq 0$, then $\omega$ has no
trivial sub-bundles.  Alternatively, $\omega$ does not admit an everywhere 
non-zero cross section.

It follows from the construction of the $p_i = \theta_1 \oplus \tau_2^{i*}(\tilde{\omega_i})$
that $\tilde{\omega_i}$ is a vector bundle over $Y_i$ with non-zero Euler class (\cite{V2}).

It will suffice to find an element $f$ of $B_i$ such that
$||\alpha(f)-f|| \geq 1$ and
\begin{displaymath}
||\mathrm{Ad}(u) \circ \alpha \circ \gamma_{im}(f)-\gamma_{im}(f)|| \geq 1
\end{displaymath}
for all unitaries $u \in B_{m}$ and natural numbers $m \in \mathbb{N}$.

Let $\tilde{f}$ be continuous function on $\mathrm{I}$ taking values in
$[0,1]$ such that $\tilde{f}(0)=0$ and $\tilde{f}(1)=1$.
Pull this function back to a function on $X_i = \mathrm{I} \times Y_i$ 
via the co-ordinate projection onto $\mathrm{I}$, keeping
the same notation.  Put $f = \tilde{f} \theta_1 \in B_i$.
Thus chosen, the element $f \in B_i$ has the desired property:
\[
||\alpha(f)-f|| \geq 1.
\]
Notice that $\theta_1 \gamma_{im}(f)
\theta_1 = (\tilde{f} \circ \pi_{im}) \theta_1$ inside $B_m$ for all natural numbers $m \geq i$,
and that $\alpha|_{B_i}(\theta_1) = \theta_1$ for every $i \in \mathbb{N}$.

Let $u$ be a unitary
element in $B_m$.  We claim that there is a 
$y_0 \in Y_m$
such that conjugation by $u$ fixes the corner
\begin{displaymath}
\theta_1(\mathrm{C}(X_m) \otimes \mathcal{K})\theta_1
\end{displaymath}
of $B_m$ at $(0,y_0) \in X_m = \mathrm{I} \times Y_m$, i.e.,
\begin{displaymath}
(u^* \theta_1 g \theta_1 u)(0,y_0) = (\theta_1 g \theta_1)(0,y_0)
\end{displaymath}
for all $g \in \mathrm{C}(X_m \otimes \mathcal{K})$.
Let $\Gamma = (x,y) \mapsto v_{(x,y)}$ be an everywhere non-zero cross 
section of $\theta_1$ over
$\{0\} \times Y_m \subseteq X_m$.  Suppose that there is no point 
$(0,y_0)$ as above.
Let $R_{(x,y)}$ denote the fibre of the
vector bundle corresponding to ${p_m|}_{\{0\} \times Y_m}$ at 
$(0,y)$, and
let $W_{(x,y)}$ denote the subspace of $R_{(x,y)}$
corresponding to $\tilde{\omega}_m$.  By assumption, the angle between $v_{(x,y)}$
and $u^* v_{(x,y)}$ is non-zero for every $(0,y) \in \{0\} \times Y_m$.  
But this implies that the projection of 
$u^* v_{(x,y)}$ onto $W_{(x,y)}$ is an everywhere non-zero
cross section of $\tilde{\omega}_{i+1}$, contradicting 
\mbox{$e(\tilde{\omega}_{i+1}) \neq 0$} and proving the claim.

Let $(0,y_0)$ be a point in $\{0\} \times Y_m$ at which $u$ fixes the 
corner
\begin{displaymath}
\theta_1(\mathrm{C}(X_m) \otimes \mathcal{K})\theta_1.
\end{displaymath}
Then,
\begin{displaymath}
(\mathrm{Ad}(u) \circ \alpha \circ \gamma_{im}(f))(0,y_0) =
\theta_1\alpha \circ \gamma_{im}(f)(0,y_0)\theta_1 \oplus g(0,y_0),
\end{displaymath}
where $g \in \omega_m B_m \omega_m$.
We conclude that
\begin{displaymath}
||\gamma_{im}(f)-\mathrm{Ad}u \circ \alpha \circ \gamma_{im}(f)||
\end{displaymath}
is bounded below by
\begin{displaymath}
\begin{array}{rl}
&||\tilde{f}(\pi_{im}(0,y_0))\theta_1  - \alpha(\tilde{f}(\pi_{im}(0,y_0)\theta_1)|| \\
= & ||\tilde{f}(0,y^{'}) - \tilde{f}(\tilde{\psi}(0,y^{'}))|| \\ 
= & 1,
\end{array}
\end{displaymath}
as desired.

Note that $\psi_i$ is homotopic to the identity
map on $B_i$ via unital endomorphisms of $B_i$ for all $i \in \mathbb{N}$ --- it is the composition two maps:  
the first is an automorphism of $B_m$ induced by a map on $X_m$, which is itself 
homotopic to the identity map on $X_m$;  the second is an inner automorphism implemented
by a unitary in the connected component of $1 \in B_m$.  Thus, $\alpha$ induces the
identity map on any $F \in \mathcal{F}_{\mathbf{R}}$ --- the restriction to functors
whose target category consists of $R$-modules is sufficient to ensure that an inductive 
limit morphism in the target category uniquely determines an automorphism of a fixed limit object.  
Since $B$ has a unique trace, $\alpha$ also induces the identity map on $\mathrm{Ell}(B)$.

Following \cite{NT}, one sees that
the absence of topological $\mathrm{K}_1$ and the fact that $\alpha$ induces 
the identity map on the Elliott invariant force $\alpha$ to induce the identity map at 
the level of the Hausdorffized algebraic $\mathrm{K}_1$-group.

The $\mathrm{KK}$-class
of $\alpha$ is the same as that of the identity map on $B$ by virtue of its 
inducing the identity map on topological $\mathrm{K}$-theory --- 
since $B$ is in the bootstrap class, $\mathrm{K}_0 B$ is free, and $\mathrm{K}_1 B = 0$
we have that
\begin{displaymath}
\mathrm{KK}^{*}(B,B) \simeq \mathrm{Hom}(\mathrm{K}_* B, \mathrm{K}_* B)
\end{displaymath}
by the Universal Coefficient Theorem (\cite{RS}).   

The stable and the real rank of a C$^*$-algebra are not relevant to the problem of
distinguising automorphsims of the algebra.  The automorphism $\alpha$ squares to 
the identity map on $B$, whence the various notions of entropy for automorphisms
of C$^*$-algebras cannot distinguish it from the identity map.  
\end{proof}

It is not clear to the author
whether the Cuntz semigroup can distinguish $\alpha$ from the identity map on $B$,
although it seems plausible.  One can, with some industry, modify the construction
of $B$ so that there exists an embedding $\iota:S_{\infty} \to \mathrm{Aut}(B)$
with the following properties:  the induced map 
\[
\overline{\iota}:S_{\infty} \to 
\mathrm{Out}(B) := \mathrm{Aut}(B)/\overline{\mathrm{Inn}(B)}
\]
is a monomorphism,
and, for each $g \in S_{\infty}$, $\iota(g)$ acts trivially on each $F \in \mathcal{F}_{\mathbf{R}}$.  
The information which goes undetected by $\mathcal{F}_{\mathbf{R}}$ is thus complicated indeed.

\section{Some remarks on the classification problem}

A classification theorem for a category $\mathcal{C}$ amounts
to proving that $\mathcal{C}$ is equivalent to a second concrete
category $\mathcal{D}$ whose objects and morphisms are well understood.
Take, for instance, the case of AF algebras:  the category $\mathcal{C}$
has AF algebras as its objects and approximate unitary equivalence classes
of isomorphisms as its morphisms, while the equivalent (classifying) category
$\mathcal{D}$ has dimension groups as its objects and order isomorphisms
of such as its morphisms.  If one does not understand $\mathcal{D}$ any 
better than $\mathcal{C}$, then one has achieved little;  the range 
of a classifying invariant is an essential part of any classification
result.

Theorems 1.1 and 1.2 show that any classifying invariant for simple
nuclear separable C$^*$-algebras will either be discontinuous with
respect to inductive limits, or not homotopy invariant even modulo
traces.  A discontinuous classifying invariant would all but exclude the
possibility of obtaining its range;  existing range results for $\mathrm{Ell}(\bullet)$
require its continuity.  The only current candidates for
non-homotopy invariant functors from the category of C$^*$-algebras
which are not captured by $\mathcal{F}$ are
the Cuntz semigroup $W(\bullet)$ or its Grothendieck enveloping group.  
Neither of these invariants is continuous with respect
to inductive limits, but this defect can perhaps be repaired by considering
these invariants as objects in the correct category.  An invariant obtained in this manner
would, while exceedingly fine, have at least the advantage of continuity with
respect to countable inductive limits.  On the other hand, the question of 
range for such an invariant remains daunting, as the following lemma shows.

\begin{lms} Let $S^{n_1},\ldots,S^{n_k}$ be a finite collection of spheres. 
Put 
\[
Y = S^{n_1} \times \cdots \times S^{n_k}, \ N = k + \sum_{i=1}^k n_i, 
\]
and let $D(Y)$ be the semigroup of Murray-von Neumann equivalence classes 
of projections in $\mathrm{M}_{\infty}(\mathrm{C}(Y))$.  Then, there is an 
order embedding 
\[
\iota:D(Y) \to W\left(\mathrm{C}\left([0,1]^N\right)\right).
\]
\end{lms} 
 
\begin{proof}
$S^{n_i}$ can be embedded more or less canonically into $[0,1]^{n_i+1}$
as the $n_i$-sphere with centre $\left(\frac{1}{2}, \ \frac{1}{2},\ldots, \frac{1}{2}\right)$ and radius $\frac{1}{4}$.
Let $S_0^{n_i} \subseteq [0,1]^{n_i+1}$ be the hollow ball
\[
S_0^{n_i} \stackrel{\mathrm{def}}{=} \left\{ \overline{x} \in [0,1]^{n_i+1} \ : \
\frac{1}{8} < \mathrm{dist}\left(\overline{x},\left(\frac{1}{2}, \ \frac{1}{2},\ldots, \frac{1}{2}\right)\right) < \frac{3}{8} \right\},
\]
and let
\[
\pi_i: S_0^{n_i} \to S^{n_i}
\]
be the projection along rays emanating from $\left(\frac{1}{2}, \ \frac{1}{2},\ldots, \frac{1}{2}\right) \in
[0,1]^{n_i+1}$.  Put
\[
Y_0 = S_0^{n_1} \times \cdots \times S_0^{n_k} \subseteq [0,1]^N; \ \
\pi = \pi_1 \times \cdots \times \pi_k.
\]
Notice that for every natural number $n$, $\mathrm{M}_n \otimes \mathrm{C}_0(Y_0)$ is
a hereditary subalgebra of $\mathrm{M}_n \otimes \mathrm{C}\left( [0,1]^N \right)$.  Let $p,q \in
\mathrm{M}_n \otimes \mathrm{C}(Y)$ be projections, and let $\pi^*(p), \pi^*(q)$ be their pullbacks
to $Y_0$.  Let $f \in \mathrm{M}_n \otimes \mathrm{C}\left([0,1]^N\right)$ be a scalar 
function taking values in $[0,1]$ which vanishes off $Y_0$ and is 
equal to one on $Y$.  Then, $f\pi^*(p), f\pi^*(q)$  are positive
elements of $\mathrm{C}\left([0,1]^N\right)$.  If $f\pi^*(p)$ and $f\pi^*(q)$ are Cuntz equivalent, then upon restriction to $Y$ we
have that $p$ and $q$ are Cuntz equivalent.  This in turn implies that $p$ and $q$ are Murray-von 
Neumann equivalent.  Now suppose that $p$ and $q$ are Murray-von Neumann equivalent.  Since this 
implies Cuntz equivalence, there exist sequences $(x_i)$ and $(y_i)$ in $\mathrm{M}_n \otimes \mathrm{C}(Y)$
such that
\[
x_i p x_i^* \stackrel{i \to \infty}{\longrightarrow} q; \ \ y_i q y_i^* \stackrel{i \to \infty}{\longrightarrow} p.
\]
Let $(g_i)$ be an approximate unit of scalar functions for $\mathrm{M}_n \otimes \mathrm{C}_0(Y_0)$.
It follows that
\[
g_i \pi^*(x_i) f \pi^*(p) \pi^*(x_i^*) g_i \stackrel{i \to \infty}{\longrightarrow} f \pi^*(q)
\]
and
\[
g_i \pi^*(y_i) f \pi^*(q) \pi^*(y_i^*) g_i \stackrel{i \to \infty}{\longrightarrow} f \pi^*(p),
\]
whence $\pi^*(p)$ and $\pi^*(q)$ are Cuntz equivalent.  The desired embedding is
\[
\iota([p]) \stackrel{\mathrm{def}}{=} \langle f \pi^*(p) \rangle.
\]
\end{proof}

\noindent
Lemma 5.1 shows that the problem of determining $W\left(\mathrm{C}\left([0,1]^N\right)\right)$
for general $N \in \mathbb{N}$ is at least as difficult as determining the isomorphism classes 
of all complex vector bundles over an arbitrary Cartesian product of spheres; this, in turn, is
a difficult unsolved problem in its own right.
Any attempt to use $W(\bullet)$
to prove a classification theorem for, say, all simple, unital AH algebras 
--- even, as Theorem 1.1 shows, if one restricts to limits of full matrix
algebras over contractible spaces, a class for which the ranges of $\mathrm{Ell}(\bullet)$,
$\mathrm{sr}(\bullet)$, $\mathrm{rr}(\bullet)$, \mbox{$\mathrm{K}$-theory} with coefficients,
and the Hausdorffized algebraic $\mathrm{K}_1$-group are known --- 
will not enjoy a salient advantage of the slow dimension growth case:
the luxury of building blocks whose invariants can be easily
and concretely described.  (Other technical obstacles are also sure to be 
much more complicated than those faced in
the work of Elliott, Gong, and Li, and their proof already runs to several 
hundred pages.)  The Cuntz semigroup is at once necessary for classification,
and unlikely to admit a range result.  

But rather than end on a pessimistic note, we enjoin the reader to view our 
results as further evidence that the Elliott invariant \emph{will} turn out to be
complete for a sufficiently well behaved class of C$^*$-algebras.  We have proved that the moment one
relaxes the slow dimension growth condition for AH algebras (and therefore,
\emph{a fortiori} for ASH algebras), one obtains counterexamples to
the Elliott conjecture of a particularly forceful nature, so that 
slow dimension growth is connected essentially to the classification problem.  
There is evidence  that slow dimension growth and 
$\mathcal{Z}$-stability are equivalent for ASH algebras --- in the case
of simple and unital AH algebras with unique trace this has recently been proved (\cite{TW1}, \cite{TW2}).  
Optimistically, $\mathcal{Z}$-stability is an abstraction of slow dimension growth, and 
the Elliott conjecture will be confirmed for all simple, separable, and nuclear
C$^*$-algebras having this property.

\end{document}